\documentclass[12pt,a4paper,twoside]{article}

\newcommand{\D}[2]{\frac{\partial #2}{\partial #1}}
\newcommand{\DD}[2]{\frac{\partial^2 #2}{\partial #1^2}}

\newcommand{\DDD}[2]{\frac{\partial^3 #2}{\partial #1^3}}
\newcommand{\DDDD}[2]{\frac{\partial^4 #2}{\partial #1^4}}

\renewcommand{\vec}[1]{\mbox{\boldmath$#1$}}
\newcommand{\Ord}[1]{{\cal O}\left(#1\right)}

\usepackage{url}
\usepackage{graphicx}

\newcommand{\pde}{{\textsc{pde}}}
\newcommand{\cM}{{\cal M}}
\newcommand{\ur}{u^+}
\newcommand{\ul}{u^-}
\newcommand{\up}{u_{j+1}}
\newcommand{\um}{u_{j-1}}
\newcommand{\umm}{u_{j-2}}
\newcommand{\upp}{u_{j+2}}
\newcommand{\uj}{u_{j}}
\newcommand{\xjmh}{\xi_{j-1/2}}
\newcommand{\xjph}{\xi_{j+1/2}}

\newcommand{\cA}{{\cal A}}

\newcommand{\KS}{Ku\-ra\-mo\-to-Siva\-shin\-sky}
\newcommand{\reduce}{\textsc{reduce}}

\begin{document}

\title{\sf Holistic finite differences accurately model the dynamics of the
Kuramoto-Sivashinsky equation }
\author{T. MacKenzie   \and
A.J. Roberts \thanks{Dept Maths \& Comp, University of Southern
Qld, Toowoomba, Qld 4350, Australia.
\protect\url{mailto:mackenzi@usq.edu.au} and
\protect\url{mailto:aroberts@usq.edu.au} respectively.}}
\maketitle

\begin{abstract}
We analyse the nonlinear Kuramoto-Sivashinsky equation to develop
an accurate finite difference approximation to its dynamics. The
analysis is based upon centre manifold theory so we are assured
that the finite difference model accurately models the dynamics
and may be constructed systematically. The theory is applied after
dividing the physical domain into small elements by introducing
insulating internal boundaries which are later removed. The
Kuramoto-Sivashinsky equation is used as an example to show how
holistic finite differences may be applied to fourth order,
nonlinear, spatio-temporal dynamical systems. This novel centre
manifold approach is holistic in the sense that it treats the
dynamical equations as a whole, not just as the sum of separate
terms.
\end{abstract}


\section{Introduction}

We continue exploring the new approach to finite difference
approximation introduced by Roberts \cite{Roberts98a} by
approximating the dynamics of solutions to the \KS{}
equation \cite{Kuramoto78,Hyman86b,Pomeau84}. In some
non-dimensional form we take the following partial differential
equation (\pde) to govern the evolution of $u(x,t)$:
\begin{equation}
    \D tu+u\D xu+R\DD xu+\DDDD xu=0\,.
    \label{Epde}
\end{equation}
This model equation includes the mechanisms of linear growth
$u_{xx}$ controlled by the parameter $R$, high-order dissipation,
$u_{xxxx}$,
and nonlinear advection/steepening, $uu_x$.  Consider implementing
the method of lines by discretising in $x$ and integrating in time
as a set of ordinary differential equations.  A finite difference
approximation to~(\ref{Epde}) on a regular grid in $x$ is
straightforward, say $x_j=jh$ for some grid spacing $h$.  For
example, the linear term
\begin{displaymath}
    \DD xu = \frac{\up-2\uj+\um}{h^2}+\Ord{h^2}\,.
\end{displaymath}
However, there are differing valid alternatives for the nonlinear
term $uu_x$: two possibilities are
\begin{equation}
    u\D xu= \frac{\uj(\up-\um)}{2h}+\Ord{h^2}
    =\frac{\up^2-\um^2}{4h}+\Ord{h^2}\,.
    \label{Ealt}
\end{equation}
Which is better?  The answer depends upon how the discretisation
of the nonlinearity interacts with the dynamics of other terms.
The conventional approach of considering the discretisation of
each term separately does not tell us.  Instead, in order to find
the best discretisation we consider the influence of all
terms in the equation in a holistic approach.

As introduced in \cite{Roberts98a} and discussed in \S\ref{Scm},
centre manifold theory \cite[e.g.]{Carr81,Roberts97a} has
appropriate characteristics to do this.  It addresses the
evolution of a dynamical system in a neighbourhood of a marginally
stable fixed point; based upon the linear dynamics the theory
guarantees that an accurate low-dimensional description of the
nonlinear dynamics may be deduced. For example the analysis herein
supports the first approximation in~(\ref{Ealt}) but with higher-order and
nonlinear corrections.  The analysis of the \KS{}
equation~(\ref{Epde}) in \S~\ref{Sbe} favours the
discretisation
\begin{eqnarray}
    &&
    \frac{d\uj}{dt} +\left[
    \frac{\uj(-\upp+9\up-9\um+\umm)}{16h}\right.
    \nonumber \\ &&
    \quad{}+\left.\left(\frac{\up^2-\um^2}{16h}\right)-\left(\frac{
    \upp\up-\umm\um}{48h}\right)
    \right]
    \nonumber \\ &&
    \quad{}+R\left(\frac{-\upp+16\up-30\uj+16\um-\umm}{12h^2}\right)
    \nonumber \\ &&
    \quad{}+\frac{\upp-4\up+6\uj-4\um+\umm}{h^4} = 0\,,
    \label{Efd}
\end{eqnarray}
as a low-order approximation.  Provided the initial conditions are
not too extreme, centre manifold theory assures us that such a
discretisation models the dynamics of~(\ref{Epde}) to errors
$\Ord{\|u\|^3,h^2}$.  Such accuracy on a relatively coarse grid is
extremely useful for such stiff \pde{}s.  Further, because the
centre manifold is composed of actual solutions to the dynamical
system, we are assured that equation~(\ref{Efd}) models the whole
of the \KS{} equation, \emph{independent} of its algebraic form.

The discretisation~(\ref{Efd}) is a low-order approximation,
centre manifold theory also provides systematic corrections.  Analysis
to higher orders in the nonlinearity, discussed in
\S~\ref{Sbe}, shows higher order corrections to the
discretisation of the nonlinear terms.  The specific finite
difference models presented here were derived by a computer
algebra program.  Computer algebra is an effective tool because of
the systematic nature of centre manifold theory \cite{Roberts96a}.

In this preliminary exploration of the approximation of the
\KS{} equation~(\ref{Epde}) we only consider an
infinite domain or strictly periodic solutions in finite domains.
Then all elements of the discretisation are identical by symmetry and
the analysis of all elements is simultaneous.  However, if physical
boundaries to the domain of the \pde{} are present, then those
elements near a physical boundary will need special treatment.
Further research is needed on this and other issues.

\section{Centre manifold theory underpins the fidelity}
\label{Scm}

Here we describe in detail one way to place the discretisation of the
\KS{} equation~(\ref{Epde}) within the purview of
centre manifold theory.

The discretisation is established via an equi-spaced grid of
collocation points, $x_j=jh$ say, for some small spacing $h$. Here
we scale the \KS{} to the scale of the grid spacing $h$. Thus our
view of the dynamics shrinks with $h$.  This is different to the
analysis in \cite{Roberts98a} and allows the linear dynamics to be
dominated by simply the highest order spatial derivative term.  We
work on the scale of the grid by transforming~(\ref{Epde}) to the
following space and time scales: $\xi=x/h$ and $\tau=t/h^4$,
giving
\begin{equation}
     \D \tau u+h^3u\D\xi u+h^2R\DD \xi u+\DDDD \xi u=0\,.
     \label{Eks}
\end{equation}
Then the crucial step: at midpoints
$\xi_{j+1/2}=(\xi_j+\xi_{j+1})/2$ artificial boundaries are
introduced:
\begin{eqnarray}
    \left[\begin{array}{c}\D \xi\ur
    -\D \xi\ul \\
    \frac{\partial^3 \ur}{\partial \xi^3}
    -\frac{\partial^3 \ul}{\partial \xi^3}
    \end{array}\right]
    &=& \left[\begin{array}{c}0\\0\end{array}\right]\,,
    \label{Ebcsd}\\
    \left(\frac{1-\gamma}{2}\right)
    \left[\begin{array}{c}\D \xi\ur
    +\D \xi\ul \\
    \frac{\partial^3 \ur}{\partial \xi^3}
    +\frac{\partial^3 \ul}{\partial \xi^3}
    \end{array}\right]
    &=&\gamma \cA
    \left[\begin{array}{c}
    \ur-\ul\\
    \frac{\partial^2 \ur}{\partial \xi^2}
    -\frac{\partial^2 \ul}{\partial \xi^2}
    \end{array}\right]\,,
    \label{Esbci}
\end{eqnarray}
where $\ur$ is just to the right of a midpoint and $\ul$ to the
left. The introduction of the near identity operator
\begin{equation}
    \cA
    =1+\frac{\partial^2_\xi}{12}-\frac{\partial^4_\xi}{720}+
    \frac{\partial^6_\xi}{30240}+\cdots
    =\frac{\partial_\xi}{2}\coth\left(\frac{\partial_\xi}{2}\right)\,,
    \label{Efudge}
\end{equation}
ensures that high-order approximations to linear terms are
obtained exactly as discussed in \cite[\S4]{Roberts98a}:  it is
remarkable that the exactly equivalent operator works for both
Burgers' equation and the \KS{} equation.  These boundaries divide
the domain into a set of elements, the $j$th element centred upon
$\xi_j$ and of width $\Delta \xi=1$.  A non-zero value of the
parameter $\gamma$ couples these elements together so that when
$\gamma=1$ the grid scaled \KS{} \pde{}~(\ref{Eks}) is effectively
restored over the whole domain by ensuring sufficient continuity
between adjacent elements.

The application of centre manifold theory is based upon a linear
picture of the dynamics.  Adjoin the dynamically trivial equations
\begin{equation}
    \D \tau \gamma=\D \tau h=0\,,
    \label{Etriv}
\end{equation}
and consider the dynamics in the extended state space
$(u(\xi),\gamma,h)$. This is a standard trick used to unfold
bifurcations \cite[\S1.5]{Carr81} or to justify long-wave
approximations \cite{Roberts88a}.  Within each element
$u=\gamma=0$ is a fixed point. Linearized about each fixed point,
that is to an error $\Ord{\|u\|^2+\gamma^2+h^2}$, the \pde{} is
\begin{displaymath}\quad
    \frac{\partial u}{\partial \tau}=\frac{\partial^4 u}{\partial \xi^4}\,,
    \quad\mbox{s.t.}\quad
    \left.\D \xi u\right|_{\xi=\pm 1/2}=
    \left.\DDD \xi u\right|_{\xi=\pm 1/2}=0\,,
\end{displaymath}
namely the hyperdiffusion equation with essentially insulating boundary
conditions.  There are thus linear eigenmodes associated with each
element:
\begin{equation}\quad
\gamma=0\,,\quad
    u\propto \left\{
    \begin{array}{ll}
        e^{\lambda_n \tau}\cos[n\pi(\xi-\xjmh)]\,, & \xjmh<\xi<\xjph\,,  \\
        0\,, & \mbox{otherwise}\,,
    \end{array}\right.
    \label{Emode}
\end{equation}
for $n=0,1,\ldots$, where the decay rate of each mode is  $
\lambda_n=-{n^4\pi^4}\,;$
together with the trivial modes $\gamma=\mbox{const}$,
$h=\mbox{const}$ and $u=0$.  In a domain with $m$ elements,
evidentally all eigenvalues are negative, $-\pi^4$ or less, except
for $m+2$ zero eigenvalues: $1$ associated with each of the $m$
elements and $2$ from the trivial~(\ref{Etriv}).  Thus, provided
the nonlinear terms in~(\ref{Eks}) are sufficiently well behaved,
the existence theorem (\cite[p281]{Carr83b}
or~\cite[p96]{Vanderbauwhede89}) guarantees that a
$m+2$~dimensional centre manifold $\cM$ exists
for~(\ref{Eks}--\ref{Etriv}).  The centre manifold $\cM$ is
parameterized by $\gamma$, $h$ and a measure of $u$ in each
element, say $\uj$: using $\vec u$ to denote the collection of
such parameters, $\cM$ is written as
\begin{equation}
    u(\xi,\tau)=v(\xi;\vec u,\gamma,h)\,.
    \label{Ecmv}
\end{equation}
In this the analysis has a very similar appearance to that of finite
elements.  The theorem also asserts that on the centre manifold the
parameters $\uj$ evolve deterministically
\begin{equation}
    \dot\uj=g_j(\vec u,\gamma,h)\,,
    \label{Ecmg}
\end{equation}
where $\dot\uj$ denotes $d\uj/d\tau$, and $g_j$ is the restriction
of~(\ref{Eks}--\ref{Etriv}) to $\cM$.  In this approach the
parameters of the description of the centre manifold may be
anything that sensibly measures the size of $u$ in each
element---we simply choose the value of $u$ at the grid points,
$\uj(\tau)=u(\xi_j,\tau)$. This provides the necessary amplitude
conditions, namely that $\uj=v(\xi_j;\vec u,\gamma,h)$.\par
The above application of the theorem establishes that in principle we
may find the dynamics~(\ref{Ecmg}) of the interacting elements of the
discretisation.  A low order approximation written in unscaled
variables is given in~(\ref{Efd}).

The next outstanding question to answer is: how can we be sure
that such a description of the interacting elements does actually
\emph{model} the dynamics of the original
system~(\ref{Eks}--\ref{Etriv})?   Here, the relevance theorem of
centre manifolds, \cite[p282]{Carr83b}
or~\cite[p128]{Vanderbauwhede89}, guarantees that all solutions
of~(\ref{Eks}--\ref{Etriv}) which remain in the neighbourhood of
the origin in $(u(\xi),\gamma,h)$ space are exponentially quickly
attracted to a solution of the $m$ finite difference
equations~(\ref{Ecmg}).  For practical purposes the rate of
attraction is estimated by the leading negative eigenvalue, here
$-\pi^4$.  Centre manifold theory also guarantees that the
stability near the origin is the same in both the model and the
original.  Thus the finite difference model will be stable if the
original dynamics are stable.  After exponentially quick
transients have died out, the finite difference
equation~(\ref{Ecmg}) on the centre manifold accurately models the
complete system~(\ref{Eks}--\ref{Etriv}).

The last piece of theoretical support tells us how to approximate
the shape of the centre manifold and the evolution thereon.
Approximation theorems such as that by Carr \& Muncaster
\cite[p283]{Carr83b} assure us that upon substituting the
ansatz~(\ref{Ecmv}--\ref{Ecmg}) into the
original~(\ref{Eks}--\ref{Etriv}) and solving to some order of
error in $\|u\|$, $\gamma$ and $h$, then $\cM$ and the evolution
thereon will be approximated to the same order.  The catch with
this application is that we need to evaluate the approximations at
$\gamma=1$ because it is only then that the artificial internal
boundaries are removed.  In some applications of such an
artificial homotopy good convergence in the parameter $\gamma$
\cite{Roberts94c} has been found.  Thus although the order of
error estimates do provide assurance, the actual error due to the
evaluation at $\gamma=1$ should be also assessed otherwise.  Here
we have crafted the interaction~(\ref{Esbci}) between elements so
that low order terms in $\gamma$ recover the exact finite
difference formula for linear terms. Note that although centre
manifold theory ``guarantees'' useful properties near the origin
in $(u(\xi),\gamma,h)$ space, because of the need to evaluate
asymptotic expressions at $\gamma=1$, we have used a weaker term
elsewhere, namely ``assures''.

\section{Numerical comparisons show the effectiveness}
\label{Sbe}

We now turn to a detailed description of the centre manifold model for
the \KS{} equation~(\ref{Epde}).

The algebraic details of the derivation of the centre manifold
model~(\ref{Ecmv}--\ref{Ecmg}) are handled by computer
algebra.\footnote{The \reduce{} computer algebra source code is
available from the authors upon request.}  In an algorithm
introduced in~\cite{Roberts96a}, the program iterates to drive to
zero the residuals of the governing differential
equation~(\ref{Eks}) and its boundary
conditions~(\ref{Ebcsd}--\ref{Efudge}).  Hence by the
Approximation theorems we assuradly construct appropriate
approximations to the centre manifold model.

The finite difference model is given by the evolution on the
centre manifold.  In order to represent the spatial fourth
derivative in the \KS{} equation we need to determine the
interactions between next-nearest neighbouring elements.  Thus the
first approximation we can consider involves quadratic terms in
$\gamma$.  After returning to $x$ and $t$ variables it is
\begin{eqnarray}
    \frac{d \uj}{d t} &=&
    -\frac{\gamma R}{h^2}\left(\up-2\uj+\um\right)
    -\frac{\gamma}{2h}\uj(\up-\um)\nonumber\\&&{}
    -\frac{\gamma^2}{h^4}\left(\upp-4\up+6\uj-4\um+\umm\right)\nonumber \\&&{}
    +\frac{\gamma^2 R}{12h^2}\left(\upp-4\up+6\uj-4\um+\umm\right)
    \nonumber \\&&{}
    +\frac{\gamma^2}{48h}\left(\upp\up+3\upp\uj-3\up-3\up\uj \right.
    \nonumber \\&&\quad{}\left.
    +3\um\uj+3\um-3\umm\uj-\umm\um\right)
    \nonumber \\&&{}
    +\frac{\gamma h^2 u_j^2}{120}\left(\up-2\uj+\um\right)\nonumber\\&&{}
    +\frac{\gamma^2
    h^2}{60480}\left[\uj^2\left(-30\upp-170\up+256\uj-170\um-30\umm\right)\right.
    \nonumber \\&& \quad{}
    +\uj\left(-126\upp\up-54\up\um-126\umm\um\right)
    \nonumber \\&& \quad{}
    +\up^2\left(10\upp-20\up+235\uj\right)
    \nonumber \\&& \quad{}
    \left.+\um^2\left(10\um-20\um+235\uj\right)\right]
    \nonumber \\&&{}
    +\Ord{\|u\|^4,\gamma^3,h^4}\,.
     \label{Ecmg24}
\end{eqnarray}
The first two lines recorded here, when evaluated for $\gamma=1$,
form the conventional second-order finite difference equation for
the \KS{} equation~(\ref{Epde}). The third line when
evaluated for $\gamma=1$ gives the fourth order accurate
corrections to the $Ru_{xx}$ term. The fourth and further lines
above start accounting systematically for the variations in the
field $u$ within each element and how they affect the evolution
through the nonlinear term.
Finite difference equations derived via this approach holistically
model all the interacting dynamics of the entire \pde.

\begin{figure}[tbp]
    \centering
    \includegraphics[width=0.9\textwidth]{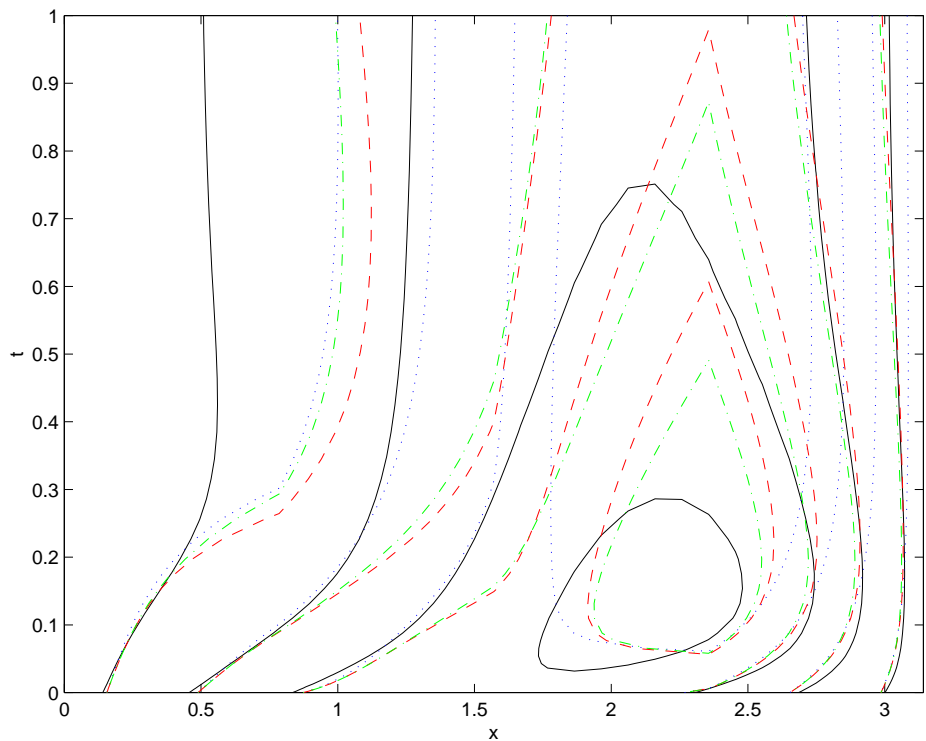}
    \caption{Contours
    of an accurate solution $u(x,t)$, -----, to compare with numerical
    approximations~(\ref{Ecmg24}): $\cdots$, the conventional
    approximation, errors ${\cal O}\left(h^2\right)$; $-~-~-$, the first
    correction, errors ${\cal O}\left(\|u\|^3\right)$; $-\cdot-\cdot-$
    , the second correction,errors ${\cal O}\left(\|u\|^4\right)$.
    \KS{} equation~(\ref{Epde}) with
    parameter $R=2$ is discretised on just $m=8$
    elements in $[0,2\pi)$ and drawn with contour interval $\Delta
    u=3$.}
    \label{F286}
\end{figure}
To show the effectiveness of the approach we compare the finite
difference model obtained from various truncations
of~(\ref{Ecmg24}) to accurate numerical solutions obtained on a
much finer grid. Choosing $m$ intervals on $[0,2\pi)$ gives an
element length $h=2\pi/m$ and grid points $x_j=jh$ for
$j=0,\ldots,m-1$. Because of the antisymmetry in $u(x,t)$ about
$x=k\pi$, when starting from the initial condition $u(x,0)=10\sin x$,
we only display the interval $[0,\pi]$. There are three different
approximations from~(\ref{Ecmg24}), with $\gamma=1$, depending
upon where the expansion is truncated.
The first two lines form a model with ${\cal O}\left(h^2\right)$
errors (a conventional finite difference approximation), the first
five lines provide the first correction with ${\cal O}
\left(\|u\|^3\right)$ errors, and all shown terms form the model
with ${\cal O}\left(\|u\|^4\right)$ errors. The solutions of these
models over $0<t<1$ with $m=8$ and $R=2$ are shown in
Figure~\ref{F286}. Observe that the leading approximation (dotted)
is significantly in error whereas the next two refinements
(dot-dashed and dashed) are overall more accurate, especially near
the peak. Such accuracy is remarkable considering the
nonlinearity, and the few points in the discretisation, $m=8$.

\section{Conclusion}

Centre manifold theory is a powerful new approach to deriving finite
difference models of dynamical systems.  Many details need to
researched for a general application of the theory.  However, there
are many promising features of this application to the \KS{}
equation~(\ref{Epde}) and the earlier example of Burgers' equation
\cite{Roberts98a}.

\bibliographystyle{plain}
\bibliography{more2}

\end{document}